\documentclass{amsart}

\usepackage{graphicx}
\usepackage[utf8]{inputenc}
\usepackage{amssymb}
\usepackage{epstopdf}
\usepackage{amsmath}
\usepackage{xcolor}
\newif\ifshowchanges
\showchangestrue

\usepackage[sort,nocompress]{cite}
\usepackage[colorlinks=true,citecolor=blue]{hyperref}
\newtheorem{theorem}{Theorem}[section]
\newtheorem{lemma}[theorem]{Lemma}

\newtheorem{definition}[theorem]{Definition}
\newtheorem{remark}[theorem]{Remark}

\numberwithin{equation}{section}
\newtheorem{prop}[theorem]{Proposition}

\def\R{{\mathbb R}}

\def\be{\begin{equation}}
\def\ee{\end{equation}}
\def\bea{\begin{eqnarray}}
\def\eea{\end{eqnarray}}
\def\beaa{\begin{eqnarray*}}
\def\eeaa{\end{eqnarray*}}

\def\ifl{\iffalse}

\def\R{{\mathbb R}}
\def\0{\Omega}

\def\d{\,{\rm d}}

\makeatletter
\AddToHook{cmd/@settitle/before}{\global\topskip=24pt\relax}
\@namedef{subjclassname@2020}{\textup{2020} Mathematics Subject Classification}
\makeatother

\begin{document}
\title[Stabilization of water waves with constant vorticity]
{Stabilization of gravity-capillary water waves with constant vorticity}

\author[Q. Ding]
{Qixing Ding}
\address{School of Mathematics and Statistics, Shanghai Polytechnic University, Shanghai 201209, China}
\email{qixingding@163.com $(\textrm{Q. Ding})$}

\author[J.-C. Zhao]
{ Jia-Cheng Zhao}
\address{School of Mathematical Sciences, Shenzhen University, Shenzhen 518061, China}
\email{zjc@szu.edu.cn $(\textrm{J.-C. Zhao})$}


\begin{abstract}
We establish an exponential stabilization result for two-dimensional periodic gravity-capillary water waves of finite depth with constant vorticity, using an external pressure localized on the free surface. The pressure is prescribed as the product of a nonnegative cut-off function and the time derivative of the surface elevation and serves as a feedback control. 
To prove the main result,
we derive the energy dissipation identity and estimate the time integral of the energy relative to the flat Couette flow, by means of the multiplier method and periodic Pohozaev identities.
A remaining lower-order term led by the vorticity is removed by a compactness argument based on unique continuation for the linearized equations. Under a condition relating the cut-off to gravity and surface tension, the relative energy decays exponentially on any interval on which a regular solution remains uniformly small in Sobolev norms. We do not assume the evenness on the solution and smallness on the vorticity.
\end{abstract}
\maketitle
\enlargethispage{2pt}

\setcounter{tocdepth}{1}
\tableofcontents

\section{Introduction}

In this paper we consider two-dimensional periodic water waves with constant vorticity  in a fluid of finite depth, bounded above by a free surface and below by a flat impermeable bed. The fluid is inviscid and incompressible with constant density, and both gravity and surface tension are taken into account. The motion is governed by the Euler equations, together with the kinematic and dynamic conditions on the free surface. More precisely, we use the Zakharov--Craig--Sulem formulation \cite{Zakharov68,CraigSulem93}, which reads
\begin{equation}\label{Waterwave}
\begin{cases}
\eta_t=G(\eta)\psi+\gamma\eta\eta_x,\\[1mm]
\psi_t+g\eta+N(\eta)\psi-\kappa H(\eta)-\gamma\eta\psi_x-\gamma\partial_x^{-1}G(\eta)\psi=-P_{\rm ext}
\end{cases}
\end{equation}
on the torus $\mathbb T=\mathbb R/(2\pi\mathbb Z)$. Here, $\gamma\in\mathbb R$ denotes the constant vorticity, $g>0$ is the gravitational acceleration, $\kappa>0$ is the surface tension coefficient, and $P_{\rm ext}$ is the external pressure acting on the free surface. The Dirichlet--Neumann operator $G(\eta)\psi$, nonlinear term $N(\eta)\psi$, and surface curvature $H(\eta)$ are formulated in Section~\ref{Section-model} below.

The mathematical theory of irrotational water waves has been extensively developed, including the existence and qualitative properties of steady traveling waves and the well-posedness of the time-dependent problem. For the Cauchy problem, see \cite{Wu97,Wu99,Lannes05,ABZ11,Wu11}. In the irrotational case, the velocity is the gradient of a harmonic potential. For two-dimensional irrotational gravity waves of infinite depth, an integrable Birkhoff normal form of order four was identified in \cite{DZ94,CW95} and rigorously justified for periodic waves in \cite{BFP23}. We refer the reader to the monograph \cite{Lannes13} for an account of the mathematical theory of water waves.

In the ocean, however, waves often interact with currents of nonzero vorticity. The study of rotational water waves is therefore important for understanding interactions between waves and currents. Constantin and Strauss \cite{CS04} used local and global bifurcation theory to establish the existence of a global continuum of steady periodic gravity waves with general vorticity. For constant vorticity, Constantin and Varvaruca \cite{CV11} studied regularity and local bifurcation, while Constantin, Strauss and Varvaruca \cite{CSV16} developed a global bifurcation theory allowing critical layers. Qualitative properties of rotational periodic waves, including their symmetry, were studied in \cite{ConEhrWah07,ConEsc04}. Related existence results for fixed-depth rotational equatorial waves were obtained in \cite{CDY20}. Rotational effects also play an important role in large-scale geophysical flows on the rotating sphere; see \cite{CG22,CJ16-1,CJ16-2}. For the time-dependent problem with constant vorticity, the velocity can be decomposed into a Couette flow and the gradient of a harmonic potential. In this formulation, Berti, Maspero and Murgante \cite{BMM24} proved almost global existence for small smooth periodic gravity--capillary waves, for almost every value of the surface tension coefficient. Their proof uses a Birkhoff normal form reduction that preserves the Hamiltonian structure of the unforced equations.

Besides the existence and long-time behavior of rotational water waves, their stability is an important question. Benjamin and Feir \cite{BF67} revealed the modulational instability of periodic gravity waves in deep water. For steady waves, variational formulations based on energy, mass, momentum and vorticity invariants were developed in \cite{CSS06} and extended to two-layer flows in \cite{CDE21}. Constantin and Strauss \cite{CS07} obtained linear stability results under restrictions on the vorticity and boundary perturbations, and a formal stability result for small-amplitude long waves with small vorticity. In the equatorial model of \cite{CDY20}, formal stability was proved for the laminar flow from which periodic waves bifurcate. Groves and Wahl\'en \cite{GW15} established conditional energetic stability of the set of constrained energy minimizers for solitary gravity--capillary waves with constant vorticity. For homogeneous water waves, the total energy consists of kinetic energy, gravitational potential energy and surface energy. Under constant atmospheric pressure, this energy is conserved. The Hamiltonian formulation of irrotational water waves was introduced by Zakharov \cite{Zakharov68}. For constant vorticity, Constantin, Ivanov and Prodanov \cite{CIP08} derived a nearly-Hamiltonian formulation, while Wahl\'en \cite{Wahlen07} obtained a canonical Hamiltonian formulation. Hamiltonian formulations were also developed for interactions between waves and currents in stratified rotational flows \cite{CIM}. When an external pressure is applied on the free surface, however, the energy changes through the work done by the pressure.

One may therefore ask whether a pressure depending on the evolving wave can make the energy decay. This is the stabilization problem, which arises in the absorption of waves in numerical and experimental wave tanks \cite{AlazardGravity18}. 
It differs from exact controllability, which asks whether prescribed states can be connected in finite time. General accounts of these questions are given in \cite{Coron07,Komornik94,Lions88}.
For incompressible fluids,
boundary controllability for the two-dimensional Euler equations was established in \cite{Coron96}, and local exact controllability to trajectories for the Navier--Stokes system with localized control was proved in \cite{FCGIP04}; further results for viscous flows can be found in \cite{CMS20,CMSZ19,LSZ22}. Control and stabilization have also been studied for dispersive equations; see, e.g.  \cite{Rosier97,LRZ10} for KdV equations and \cite{Lebeau92,Laurent10,Laurent10-1} for Schr\"odinger equations. 

Water-wave problems bring these two settings together: the fluid motion is governed by the incompressible Euler equations, while the evolution of the free surface has a dispersive feature. The moving boundary and the resulting nonlocal, quasilinear surface equations introduce additional difficulties for control and stabilization. Boundary control of linear water waves was studied by Reid and Russell \cite{RR85}. The related problem of generating prescribed waves by a localized pressure was studied by Alazard, Baldi and Han-Kwan \cite{ABH18}, who established local exact controllability of  two-dimensional irrotational gravity--capillary waves. Zhu \cite{Zhu20} extended this result to higher dimensions under a geometric control condition. More recently, Haus and Pasquali \cite{HP26} proved local exact controllability for two-dimensional irrotational gravity--capillary waves over an open and dense class of admissible bottom topographies, and for waves with arbitrary constant vorticity over a flat bottom. For irrotational gravity waves without surface tension, Alazard \cite{AlazardGravity18} obtained energy decay estimates using a feedback prescribed through the horizontal derivative of the pressure in an absorbing zone. In the presence of surface tension, Alazard \cite{Alazard17} used the feedback $P_{\rm ext}=\chi\eta_t$ and proved an integral estimate yielding exponential energy decay for even periodic regular solutions satisfying geometric bounds and a compatibility condition between the cut-off, gravity and surface tension.

In the present work, we consider a localized pressure feedback for waves with constant vorticity, with the aim of damping perturbations of the flat Couette flow. 
For regular solutions, we prove exponential decay of the energy relative to this flow under the feedback $P_{\rm ext}=\chi\eta_t$, on every interval where $(\eta,\psi-\langle\psi\rangle)$ remains uniformly small in $H^{s+1/2}\times H^s$, with $s>5/2$.
The estimates allow any fixed vorticity $\gamma$, subject to the stated condition on the cut-off. The main result is precisely stated in Theorem~\ref{Thm-stabilization} later. 
This extends the stabilization result of \cite{Alazard17} (addressing the case of $\gamma=0$) to the gravity-capillary waves with nonzero vorticity.
Our main contributions are as follows.
\begin{itemize}
\item \emph{Stabilization near a Couette flow.} We prove exponential decay of the energy relative to the flat Couette flow under the localized pressure feedback  used in \cite{Alazard17}. The vorticity is fixed but need not be small. 
\vspace{1mm}

\item \emph{No evenness assumption.} We treat periodic solutions without the evenness assumption associated with the tank setting in \cite{Alazard17}. The multiplier and Pohozaev calculations are carried out over a full spatial period, with the additional vorticity terms explicitly estimated.

\vspace{1mm}

\item \emph{Removal of a lower-order remainder.} Estimating the nonlocal vorticity term $\gamma\partial_x^{-1}G(\eta)\psi$ produces a lower-order $L^2$ remainder whose coefficient depends on $\gamma$ and need not be small. The localized kinetic energy estimate produces a remainder of the same type. We remove the combined remainder by proposing a compactness argument based on unique continuation for the linearized equations, thereby accomplishing the decay estimate.
\end{itemize}

The plan of the paper is as follows. In Section~\ref{Section-model}, we introduce the governing equations and derive the energy balance and the estimates for the harmonic potential. In Section~\ref{Section-multiplier}, we adapt Alazard's multiplier calculation, compare the kinetic and potential energies, and apply periodic Pohozaev identities. The periodic calculation is needed because evenness of the surface variables is not preserved for nonzero vorticity. In Section~\ref{Section-decay}, we estimate the pressure and the vorticity terms. The remaining term of lower order is removed by performing a compactness--uniqueness argument for the linearized equations. The resulting integral estimate and the energy balance yield exponential decay.

\section{Governing equations and energy profiles}\label{Section-model}

\subsection{Formulation of the model}

We consider a two-dimensional inviscid and incompressible fluid of unit density occupying
\[
D_\eta(t)=\{(x,y):-\pi<x<\pi,\ -h<y<\eta(t,x)\},
\]
where $h>0$, the vertical coordinate $y$ points upwards, and the two vertical edges are identified by periodicity. We assume $\eta>-h$ and take the surface elevation to have zero spatial mean, so that $h$ is the mean depth.
Let $(u,v)$ be the velocity and $P$ the pressure relative to the constant atmospheric pressure. The motion is governed by
\begin{equation*}
\begin{aligned}
&u_t+uu_x+vu_y=-P_x,\\
&v_t+uv_x+vv_y=-P_y-g,
\end{aligned}
\end{equation*}
with $u_x+v_y=0$ in $D_\eta(t)$, where $g>0$ is the gravitational acceleration. The boundary conditions are
\begin{align*}
v&=\eta_t+u\eta_x
&&\text{on }y=\eta(t,x),\\
P&=P_{\rm ext}-\kappa H(\eta)
&&\text{on }y=\eta(t,x),\\
v&=0
&&\text{on }y=-h.
\end{align*}
Here $P_{\rm ext}$ is the external pressure, $\kappa>0$ is the surface tension coefficient, and
\[
H(\eta)=\partial_x\left(\frac{\eta_x}{\sqrt{1+\eta_x^2}}\right)
\]
is the curvature of the free surface.
We assume constant vorticity $v_x-u_y=\gamma$, with $\gamma\in\R$, and fix the circulation by
\[
\int_{-\pi}^{\pi}\bigl(u(t,x,-h)-\gamma h\bigr)\d x=0.
\]
The velocity then admits the decomposition
\begin{equation}\label{velocity}
u=\Phi_x-\gamma y,\quad v=\Phi_y,
\end{equation}
where $\Phi$ is a $2\pi$-periodic harmonic function in $D_\eta(t)$ satisfying $\Phi_y=0$ on the bed. We introduce its surface trace and the Dirichlet--Neumann operator by
\[
\psi(t,x)=\Phi(t,x,\eta(t,x)),\quad G(\eta)\psi=\left(\Phi_y-\eta_x\Phi_x\right)\big|_{y=\eta(t,x)}.
\]
The quadratic term in the surface equations is
\[
N(\eta)\psi=\frac12\psi_x^2-\frac{\bigl(G(\eta)\psi+\eta_x\psi_x\bigr)^2}{2(1+\eta_x^2)},
\]
which also has the expression
\begin{equation}\label{Nbulk}
N(\eta)\psi=\left(\frac12\Phi_x^2-\frac12\Phi_y^2+\eta_x\Phi_x\Phi_y\right)\bigg|_{y=\eta}.
\end{equation}
We denote by $\partial_x^{-1}$ the zero-mean periodic primitive, with Fourier symbol $(\mathrm i j)^{-1}$ for $j\ne0$.
It applies to $G(\eta)\psi$ since
\begin{equation*}
\int_{-\pi}^{\pi}G(\eta)\psi\d x=0.
\end{equation*}
With a suitable choice of the time-dependent additive constant in $\Phi$, the Euler equations and boundary conditions
take the surface form, which is exactly \eqref{Waterwave}. 
When $P_{\rm ext}=0$, this is the system in \cite{BMM24}. For the surface formulation, see also \cite[Section~3]{CIP08} and \cite[Section~1.2]{HP26}.
The reference state is the flat Couette flow,
\[
\eta=0,\quad (u,v)=(-\gamma y,0).
\]
For $r\in\mathbb R$, let $H^r=H^r(\mathbb T)$ denote the Sobolev space of $2\pi$-periodic functions, and let $H^r_0$ be its zero-mean subspace. For the boundary calculations below, we use the regularity in \cite[Definition~2.1]{Alazard17}.

\begin{definition}[Regular solution]\label{Def-regularsol}
A real-valued solution $(\eta,\psi,P_{\rm ext})$ of \eqref{Waterwave} on $[0,T]$ is called a regular solution if, for some $s>5/2$,
\begin{equation*}
\eta\in C([0,T];H^{s+1/2}_0),\quad\psi\in C([0,T];H^s),\quad P_{\rm ext}\in C([0,T];H^{s-1}),
\end{equation*}
and $\eta>-h$ on $[0,T]$.
\end{definition}

The condition $s>5/2$ ensures, by Sobolev embedding, that the surface variables have continuous spatial derivatives up to order two. The surface equations also give
\[
\eta_t\in C([0,T];H^{s-1}),\quad\psi_t\in C([0,T];H^{s-3/2}).
\]
These assumptions provide the regularity used in the boundary calculations and integrations by parts below.

We now choose the external pressure so that it acts only near the ends of one period. Following \cite[Definition~2.3]{Alazard17}, fix $0<\delta<\pi$ and choose a smooth, even, $2\pi$-periodic function $\varphi$ such that $0\le\varphi\le1$, $x\varphi_x(x)\le0$ on $[-\pi,\pi]$, and
\[
\varphi(x)=1\quad\text{for }|x|\le\pi-\delta,\quad\varphi(x)=0\quad\text{for }\pi-\delta/2\le|x|\le\pi.
\]
The cut-off function $\chi$ used for the localized control is
\begin{equation}\label{cutoff}
m(x)=x\varphi(x),\quad\chi(x)=1-m_x(x).
\end{equation}
Since $m$ vanishes near the endpoints, it extends to a smooth periodic function. Moreover,
\[
\chi\ge0,\quad\int_{-\pi}^{\pi}\chi(x)\d x=2\pi,\quad\chi=0\quad\text{on }[-\pi+\delta,\pi-\delta].
\]
In the present problem, the absorbing zone is the part of the free surface parametrized by the horizontal set
\begin{equation*}
\omega_{\rm abs}=[-\pi,-\pi+\delta]\cup[\pi-\delta,\pi].
\end{equation*}
The two end intervals form a single zone under the periodic identification. The cut-off satisfies $\operatorname{supp}\chi\subset\omega_{\rm abs}$ and $\chi=1$ for $\pi-\delta/2\le|x|\le\pi$. 
As in \cite{Alazard17}, we assume
\begin{equation}\label{capillary-cutoff-condition}
\kappa\|m_{xx}\|_{L^\infty}^2\le g.
\end{equation}
This condition will be used to estimate a term involving surface tension.
We prescribe the pressure by
\begin{equation}\label{Feedback-law}
P_{\rm ext}=\chi(x)\eta_t
=\chi(x)\bigl(G(\eta)\psi+\gamma\eta\eta_x\bigr).
\end{equation}
Thus the pressure acts only on the free surface portion
$\{(x,\eta(t,x)):x\in\omega_{\rm abs}\}$.

\subsection{The energy relative to  flat Couette flow}

Over one period, the total energy is the sum of the kinetic energy, the gravitational potential energy and the surface energy:
\[
\int_{-\pi}^{\pi}\int_{-h}^{\eta(t,x)} \left(\frac12(u^2+v^2)+gy\right)\d y\d x+\kappa\int_{-\pi}^{\pi}\sqrt{1+\eta_x^2}\d x.
\]
Even when the free surface is flat and $\Phi$ is constant, the fluid has the velocity $(-\gamma y,0)$ and hence has nonzero kinetic energy if $\gamma\ne0$. We therefore measure the energy relative to this flat Couette flow. The subtracted energy depends only on the fixed parameters and does not affect the energy balance.
To express this energy in terms of $(\eta,\psi)$, we use \eqref{velocity} and write
\[
\frac12(u^2+v^2)=\frac12|\nabla\Phi|^2-\gamma y\Phi_x+\frac{\gamma^2}{2}y^2.
\]
The first term is the kinetic energy of the irrotational part. Since $\Phi$ is harmonic and satisfies the bed condition, Green's formula gives
\begin{equation}\label{kinetic}
K(t):=\frac12\int_{D_\eta(t)}|\nabla\Phi|^2\d x\d y=\frac12\int_{-\pi}^{\pi}\psi G(\eta)\psi\d x.
\end{equation}
The contributions from the vertical edges cancel by periodicity. For the cross term, integration in the periodic variable gives
\[
-\gamma\int_{-\pi}^{\pi}\int_{-h}^{\eta(t,x)}y\Phi_x\d y\d x=\gamma\int_{-\pi}^{\pi}\psi\eta\eta_x\d x=-\frac{\gamma}{2}\int_{-\pi}^{\pi}\eta^2\psi_x\d x.
\]
The change in the kinetic energy of the shear flow is
\[
\frac{\gamma^2}{2}\int_{-\pi}^{\pi}\left(\int_{-h}^{\eta(t,x)}y^2\d y-\int_{-h}^{0}y^2\d y\right)\d x=\frac{\gamma^2}{6}\int_{-\pi}^{\pi}\eta^3\d x.
\]
Similarly, subtracting the gravitational potential energy of the flat flow gives $g\int_{-\pi}^{\pi}\eta^2\d x/2$, while subtracting the energy of the flat surface gives $\kappa\int_{-\pi}^{\pi}(\sqrt{1+\eta_x^2}-1)\d x$.
Combining these expressions, we obtain
\begin{equation}\label{energy}
\begin{aligned}
\mathcal E_\gamma(\eta,\psi)={}&\frac12\int_{-\pi}^{\pi}\psi G(\eta)\psi\d x+\frac g2\int_{-\pi}^{\pi}\eta^2\d x+\kappa\int_{-\pi}^{\pi}\left(\sqrt{1+\eta_x^2}-1\right)\d x\\
&-\frac{\gamma}{2}\int_{-\pi}^{\pi}\eta^2\psi_x\d x+\frac{\gamma^2}{6}\int_{-\pi}^{\pi}\eta^3\d x.
\end{aligned}
\end{equation}
This agrees with the Hamiltonian in \cite{BMM24}, after subtracting its value at the flat Couette flow. In particular, $\mathcal E_\gamma(0,0)=0$. 
The first three terms in \eqref{energy} constitute $\mathcal E_0$.
The remaining two terms arise from the underlying shear flow.

The work done by the external pressure gives the following energy balance.

\begin{prop}\label{Prop-energyidentity}
Every regular solution of \eqref{Waterwave} satisfies
\begin{equation}\label{energyidentity}
\frac{\d}{\d t}\mathcal E_\gamma(t)=-\int_{-\pi}^{\pi}P_{\rm ext}\eta_t\d x.
\end{equation}
Under the pressure law \eqref{Feedback-law}, we have
\begin{equation}\label{dissipation}
\mathcal E_\gamma(t)+\int_a^t\mathcal D(\tau)\d\tau=\mathcal E_\gamma(a),\quad\mathcal D(t)=\int_{-\pi}^{\pi}\chi(x)\eta_t(t,x)^2\d x,
\end{equation}
for $0\le a\le t\le T$.
\end{prop}

\begin{proof}
We compute the variations of \eqref{energy} with respect to the two surface variables. Since $G(\eta)$ is self-adjoint, differentiation with respect to $\psi$ gives
\begin{equation}\label{energy-psi}
\frac{\delta\mathcal E_\gamma}{\delta\psi}=G(\eta)\psi+\gamma\eta\eta_x=\eta_t.
\end{equation}

For the variation with respect to $\eta$, we first consider $K$. If the surface is varied by $\dot\eta$ while its potential trace $\psi$ is fixed, the corresponding variation of $\Phi$ satisfies $\dot\Phi=-\Phi_y\dot\eta$ on $y=\eta$.
Here the dot denotes a variation, not a time derivative. Green's formula therefore yields
\[
\dot K=\int_{-\pi}^{\pi}\left[\frac12|\nabla\Phi|^2\big|_{y=\eta}-\Phi_y\big|_{y=\eta}\,G(\eta)\psi\right]\dot\eta\d x=\int_{-\pi}^{\pi}N(\eta)\psi\,\dot\eta\d x.
\]
Differentiating the remaining terms in \eqref{energy}, we obtain
\begin{equation}\label{energy-eta}
\frac{\delta\mathcal E_\gamma}{\delta\eta}=N(\eta)\psi+g\eta-\kappa H(\eta)-\gamma\eta\psi_x+\frac{\gamma^2}{2}\eta^2.
\end{equation}
Combining \eqref{energy-eta} with the second equation of \eqref{Waterwave}, we obtain
\begin{equation}\label{psi-energy-eta}
\psi_t+\frac{\delta\mathcal E_\gamma}{\delta\eta}=\gamma\partial_x^{-1}G(\eta)\psi+\frac{\gamma^2}{2}\eta^2-P_{\rm ext}.
\end{equation}
Using the first equation and the identity
\[
\partial_x^{-1}(\eta\eta_x)=\frac12\left(\eta^2-\langle\eta^2\rangle\right),
\]
we can rewrite \eqref{psi-energy-eta} as
\[
\psi_t+\frac{\delta\mathcal E_\gamma}{\delta\eta}=\gamma\partial_x^{-1}\eta_t+\frac{\gamma^2}{2}\langle\eta^2\rangle-P_{\rm ext}.
\]
Consequently, by \eqref{energy-psi} we have
\[
\frac{\d}{\d t}\mathcal E_\gamma=\int_{-\pi}^{\pi}\left(\frac{\delta\mathcal E_\gamma}{\delta\eta}\eta_t+\frac{\delta\mathcal E_\gamma}{\delta\psi}\psi_t\right)\d x=\int_{-\pi}^{\pi}\eta_t\left(\gamma\partial_x^{-1}\eta_t+\frac{\gamma^2}{2}\langle\eta^2\rangle-P_{\rm ext}\right)\d x.
\]
Since $\langle\eta_t\rangle=0$ and $\partial_x^{-1}$ is skew-adjoint on zero-mean periodic functions,
\[
\int_{-\pi}^{\pi}\eta_t\partial_x^{-1}\eta_t\d x=0,\quad\int_{-\pi}^{\pi}\eta_t\d x=0.
\]
The first two terms therefore vanish, which proves \eqref{energyidentity}.  Substituting $P_{\rm ext}=\chi\eta_t$ and integrating in time gives \eqref{dissipation}.
\end{proof}

The energy is therefore conserved when $P_{\rm ext}=0$ and is non-increasing under \eqref{Feedback-law}. To relate it to the surface variables, we first record the required estimates for the harmonic potential.
For a periodic function $f$, write
\[
\langle f\rangle=\frac1{2\pi}\int_{-\pi}^{\pi}f(x)\d x,\quad\widetilde f=f-\langle f\rangle.
\]
In particular, we put $\widetilde\Phi=\Phi-\langle\psi\rangle$. Then $\widetilde\Phi$ is the harmonic extension of $\widetilde\psi$.
The following estimates are the mixed-boundary elliptic estimates used for water waves; see \cite{Lannes13} and the appendices of \cite{Alazard17}. Their constants are uniform for small $\|\eta\|_{W^{2,\infty}}$, with $h$ fixed.

\begin{lemma}\label{Lemma-elliptic}
There are $\delta_0>0$ and $C>0$, depending on $h$, such that, if $\|\eta\|_{W^{2,\infty}}\le\delta_0$ and $\eta\ge-h/2$, then
\begin{align}
C^{-1}\|\widetilde\psi\|_{H^{1/2}}^2&\le K\le C\|\widetilde\psi\|_{H^{1/2}}^2,\label{Knorm}\\
\|\widetilde\Phi\|_{L^2(D_\eta)}&\le C\|\widetilde\psi\|_{L^2},\label{PhiL2}\\
\|\partial_x^{-1}G(\eta)\psi\|_{L^2}&\le C\|\widetilde\psi\|_{L^2},\label{AGbound}\\
\int_{-\pi}^{\pi}\Phi_x(x,-h)^2\d x&\le CK,\label{bottomtrace}
\end{align}
where $K$ is given by \eqref{kinetic}. In addition, 
\begin{equation}\label{DNcontinuity}
\|(G(\eta)-G(0))\psi\|_{H^{-1/2}} \le C\|\eta\|_{W^{1,\infty}}\|\widetilde\psi\|_{H^{1/2}}.
\end{equation}
\end{lemma}

\begin{proof}
Flatten the domain by
\[
(x,z)\mapsto
\left(x,z+\left(1+\frac zh\right)\eta(x)\right),\quad -h<z<0.
\]
The transformed Dirichlet form is uniformly elliptic and differs from the flat one by $O(\|\eta\|_{W^{1,\infty}})$ in its coefficients. The variational characterization of the harmonic extension and the trace theorem give \eqref{Knorm}, with the constant mode removed. Comparing the two Dirichlet problems and polarizing their forms gives \eqref{DNcontinuity}.

For \eqref{PhiL2}, solve $-\Delta z=f$ in $D_\eta$, with $z=0$ at the top and $z_y=0$ at the bottom. The $H^2$ estimate for this problem is uniform for these graph domains and implies
\[
\|\partial_n z\|_{L^2(\{y=\eta\})}\le C\|f\|_{L^2(D_\eta)}.
\]
Green's formula applied to $z$ and $\widetilde\Phi$, followed by duality, gives \eqref{PhiL2}. The surface measure on the graph is uniformly comparable to $\d x$.
The estimate $G(\eta):H^1\to L^2$ and self-adjointness give $G(\eta):L^2\to H^{-1}$. Compose this with $\partial_x^{-1}:H^{-1}_0\to L^2_0$ to obtain \eqref{AGbound}.
A flat strip of width $h/4$ adjacent to the bottom lies in $D_\eta$. The nonconstant Fourier coefficients of $\Phi$ in this strip have the form $a_j\cosh(|j|(y+h))$, because of the bottom condition. Their Dirichlet energy controls $\sum_{j\ne0}j^2|a_j|^2$. Parseval's identity proves \eqref{bottomtrace}.
\end{proof}

The two vorticity terms in \eqref{energy} have no fixed sign. The following lemma shows that they are controlled by $\mathcal E_0$ when the free surface is sufficiently close to the flat one.

\begin{lemma}\label{Lemma-coercivity}
For fixed $g,\kappa,h>0$ and $\gamma\in\R$, there are $\delta_1>0$ and $C\ge1$ such that, if $\|\eta\|_{W^{2,\infty}}\le\delta_1$ and $\langle\eta\rangle=0$, then
\begin{equation}\label{coercivity}
\frac12\mathcal E_0\le\mathcal E_\gamma\le\frac32\mathcal E_0,\quad C^{-1}\bigl(\|\eta\|_{H^1}^2+\|\widetilde\psi\|_{H^{1/2}}^2\bigr)\le\mathcal E_\gamma\le C\bigl(\|\eta\|_{H^1}^2+\|\widetilde\psi\|_{H^{1/2}}^2\bigr).
\end{equation}
\end{lemma}

\begin{proof}
For small $\|\eta_x\|_{L^\infty}$, the capillary energy is comparable to $\kappa\|\eta_x\|_{L^2}^2$. Using \eqref{Knorm} therefore gives the energy equivalence for $\mathcal E_0$. Integration by parts gives
\[
-\frac\gamma2\int_{-\pi}^{\pi}\eta^2\psi_x\d x=\gamma\int_{-\pi}^{\pi}\widetilde\psi\eta\eta_x\d x.
\]
Its absolute value is bounded by $C|\gamma|\|\eta\|_{L^\infty}\mathcal E_0$. The remaining cubic term is bounded by $C\gamma^2\|\eta\|_{L^\infty}\mathcal E_0$. Decrease $\delta_1$ to obtain \eqref{coercivity}.
\end{proof}

\section{The multiplier method}\label{Section-multiplier}

To estimate the time integral of the energy, we relate the energy of the flow to the motion of the free surface above $\omega_{\rm abs}$. We adapt the multiplier calculation of \cite[Section~3]{Alazard17}, combining the surface equations with identities for the harmonic potential. The identities in this section hold for arbitrary external pressure, and the feedback law \eqref{Feedback-law} is used in Section~\ref{Section-decay}.

For convenience, we write
\[
w=G(\eta)\psi,\quad S=\sqrt{1+\eta_x^2}.
\]
All spatial boundary integrals are taken over $[-\pi,\pi]$. The identities are written at each time and integrated over a time interval at the end of the section.

\subsection{Multiplier identities for the surface equations}

We begin by multiplying the first equation in \eqref{Waterwave} by $m\psi_x$ and the second by $-m\eta_x$. After adding the resulting expressions, we integrate in space and time. Computing the same integral by integration by parts gives the following identity.

\begin{lemma}\label{Lemma-basic}
For each fixed time, set
\[
F(t)=\int_{-\pi}^{\pi}m\psi_xw\d x+\int_{-\pi}^{\pi}m\eta_xN(\eta)\psi\d x.
\]
Then every regular solution of \eqref{Waterwave} satisfies
\begin{equation}\label{basic-multiplier}
\begin{aligned}
\int_{-\pi}^{\pi}m_x\left(-\eta\psi_t-\frac g2\eta^2-\frac{\kappa}{S}\right)\d x+F(t)&=-\frac{\d}{\d t}\int_{-\pi}^{\pi}\partial_x(m\eta)\psi\d x-\int_{-\pi}^{\pi}P_{\rm ext}m\eta_x\d x\\
&\quad+\gamma\int_{-\pi}^{\pi}m\eta_x\partial_x^{-1}w\d x.
\end{aligned}
\end{equation}
\end{lemma}

\begin{proof}
Consider
\[
A=\int_0^T\int_{-\pi}^{\pi}\left(\eta_t m\psi_x-\psi_t m\eta_x\right)\d x\d t.
\]
We first integrate the term containing $\psi_x$ by parts in $x$. Since $m$ is independent of time, a further integration by parts in $t$ gives
\begin{equation}\label{A-parts}
\begin{aligned}
A&=-\int_0^T\int_{-\pi}^{\pi}\left(\partial_x(m\eta_t)\psi+m\eta_x\psi_t\right)\d x\d t\\
&=-\left[\int_{-\pi}^{\pi}\partial_x(m\eta)\psi\d x\right]_0^T+\int_0^T\int_{-\pi}^{\pi}m_x\eta\psi_t\d x\d t.
\end{aligned}
\end{equation}
Here $[f]_0^T=f(T)-f(0)$, and the spatial boundary terms vanish by periodicity.

We next compute $A$ by substituting the two equations in \eqref{Waterwave}. The local terms involving vorticity cancel:
\[
(\gamma\eta\eta_x)m\psi_x-(\gamma\eta\psi_x)m\eta_x=0.
\]
The remaining terms give
\begin{equation}\label{A-identity}
\begin{aligned}
A={}&\int_0^T F(t)\d t+\int_0^T\int_{-\pi}^{\pi}\left(P_{\rm ext}+g\eta-\kappa H(\eta)\right)m\eta_x\d x\d t\\
&-\gamma\int_0^T\int_{-\pi}^{\pi}m\eta_x\partial_x^{-1}w\d x\d t.
\end{aligned}
\end{equation}
The gravity term is integrated by parts as
\[
\int_{-\pi}^{\pi}g\eta m\eta_x\d x=-\frac g2\int_{-\pi}^{\pi}m_x\eta^2\d x.
\]
For the surface tension term, we use $H(\eta)=\partial_x(\eta_x/S)$ and $S_x=\eta_x\eta_{xx}/S$. Thus
\[
\begin{aligned}
-\int_{-\pi}^{\pi}H(\eta)m\eta_x\d x&=\int_{-\pi}^{\pi}\frac{\eta_x}{S}\left(m_x\eta_x+m\eta_{xx}\right)\d x=\int_{-\pi}^{\pi}m_x\frac{\eta_x^2}{S}\d x+\int_{-\pi}^{\pi}mS_x\d x\\
&=\int_{-\pi}^{\pi}m_x\left(\frac{\eta_x^2}{S}-S\right)\d x=-\int_{-\pi}^{\pi}\frac{m_x}{S}\d x.
\end{aligned}
\]
Substituting these formulas into \eqref{A-identity} and comparing with \eqref{A-parts}, we obtain the time integral of \eqref{basic-multiplier} over $[0,T]$. The same calculation holds on every subinterval. Since the solution is smooth, differentiation with respect to the upper endpoint yields \eqref{basic-multiplier}.
\end{proof}

The term $-\eta\psi_t$ in \eqref{basic-multiplier} can be eliminated using the second surface equation. This brings out the gravitational and surface energies. The integration by parts also produces a surface tension remainder, which we write as
\begin{equation}\label{R1}
R_1(t)=\kappa\int_{-\pi}^{\pi}\frac{m_{xx}\eta\eta_x}{S}\d x+\kappa\int_{-\pi}^{\pi}m_x\left(1-\frac{1+\eta_x^2/2}{S}\right)\d x.
\end{equation}
Using $m_x=1-\chi$, we separate the integral over the whole period from the one containing $\chi$. 

\begin{lemma}\label{Lemma-potential}
With $R_1(t)$ defined by \eqref{R1}, every regular solution of \eqref{Waterwave} satisfies
\begin{equation}\label{potential-identity}
\begin{aligned}
\int_{-\pi}^{\pi}\left(\frac g2\eta^2+\frac{3\kappa}{2}\frac{\eta_x^2}{S}\right)\d x+R_1&=\int_{-\pi}^{\pi}\chi\left(\frac g2\eta^2+\frac{3\kappa}{2}\frac{\eta_x^2}{S}\right)\d x-\frac{\d}{\d t}\int_{-\pi}^{\pi}\partial_x(m\eta)\psi\d x\\
&\quad-\int_{-\pi}^{\pi}\partial_x(m\eta)(P_{\rm ext}+N(\eta)\psi)\d x-\int_{-\pi}^{\pi}m\psi_xw\d x\\
&\quad+\gamma\int_{-\pi}^{\pi}\partial_x(m\eta)\partial_x^{-1}w\d x+\gamma\int_{-\pi}^{\pi}m_x\eta^2\psi_x\d x.
\end{aligned}
\end{equation}
\end{lemma}

\begin{proof}
We substitute the second equation of \eqref{Waterwave} for $\psi_t$ in \eqref{basic-multiplier}. After integrating the surface tension term by parts, we obtain
\begin{align*}
&\int_{-\pi}^{\pi}m_x\left(\frac g2\eta^2+\kappa\frac{\eta_x^2}{S}-\frac\kappa S\right)\d x+\kappa\int_{-\pi}^{\pi}\frac{m_{xx}\eta\eta_x}{S}\d x\\
&\quad=-\frac{\d}{\d t}\int_{-\pi}^{\pi}\partial_x(m\eta)\psi\d x-\int_{-\pi}^{\pi}\partial_x(m\eta)(P_{\rm ext}+N(\eta)\psi)\d x-\int_{-\pi}^{\pi}m\psi_xw\d x\\
&\qquad+\gamma\int_{-\pi}^{\pi}\partial_x(m\eta)\partial_x^{-1}w\d x+\gamma\int_{-\pi}^{\pi}m_x\eta^2\psi_x\d x.
\end{align*}
Since $\int_{-\pi}^{\pi}m_x\d x=0$, the left-hand side is
\[
\int_{-\pi}^{\pi}m_x\left(\frac g2\eta^2+\frac{3\kappa}{2}\frac{\eta_x^2}{S}\right)\d x+R_1.
\]
Writing $m_x=1-\chi$ now gives \eqref{potential-identity}.
\end{proof}

To introduce the kinetic energy with a positive coefficient, we derive a weighted relation by differentiating $\int_{-\pi}^{\pi}\vartheta\eta\psi\d x$ in time.

\begin{lemma}\label{Lemma-equipartition}
For every smooth periodic function $\vartheta=\vartheta(x)$,
\begin{equation}\label{equipartition}
\begin{aligned}
\int_{-\pi}^{\pi}\vartheta\left(g\eta^2+\kappa\frac{\eta_x^2}{S}\right)\d x&=\int_{-\pi}^{\pi}\vartheta\psi w\d x-\frac{\d}{\d t}\int_{-\pi}^{\pi}\vartheta\eta\psi\d x-\int_{-\pi}^{\pi}\vartheta\eta(P_{\rm ext}+N(\eta)\psi)\d x\\
&\quad-\kappa\int_{-\pi}^{\pi}\frac{\vartheta_x\eta\eta_x}{S}\d x+\gamma\int_{-\pi}^{\pi}\vartheta\eta\partial_x^{-1}w\d x+\gamma\int_{-\pi}^{\pi}\vartheta\bigl(\psi\eta\eta_x+\eta^2\psi_x\bigr)\d x.
\end{aligned}
\end{equation}
\end{lemma}

\begin{proof}
Since $\vartheta$ is independent of time, differentiating the weighted integral and using the first equation in \eqref{Waterwave} gives
\[
\begin{aligned}
\frac{\d}{\d t}\int_{-\pi}^{\pi}\vartheta\eta\psi\d x&=\int_{-\pi}^{\pi}\vartheta\psi\eta_t\d x+\int_{-\pi}^{\pi}\vartheta\eta\psi_t\d x\\
&=\int_{-\pi}^{\pi}\vartheta\psi w\d x+\int_{-\pi}^{\pi}\vartheta\eta\psi_t\d x+\gamma\int_{-\pi}^{\pi}\vartheta\psi\eta\eta_x\d x.
\end{aligned}
\]
We now substitute the second equation in \eqref{Waterwave} for $\psi_t$. Integrating the surface tension term by parts gives
\[
\int_{-\pi}^{\pi}\vartheta\eta H(\eta)\d x=-\int_{-\pi}^{\pi}\partial_x(\vartheta\eta)\frac{\eta_x}{S}\d x=
-\int_{-\pi}^{\pi}\vartheta\frac{\eta_x^2}{S}\d x-\int_{-\pi}^{\pi}\frac{\vartheta_x\eta\eta_x}{S}\d x.
\]
Inserting this identity into the expression for the time derivative and moving the gravitational and surface tension terms to the left yields \eqref{equipartition}.
\end{proof}

We now combine the potential-energy identity with the weighted relation between kinetic and potential energies. The pressure and time-derivative terms combine through the following multiplier:
\[
\zeta=\partial_x(m\eta)+\frac32(1-m_x)\eta-\frac14\eta =m\eta_x+\left(\frac54-\frac12m_x\right)\eta.
\]
The same calculation produces the energy expression
\[
\widetilde H=K+\frac g2\int_{-\pi}^{\pi}\eta^2\d x
+\kappa\int_{-\pi}^{\pi}\frac{\eta_x^2}{S}\d x.
\]
Since $S-1=\frac{\eta_x^2}{S+1}\le\frac{\eta_x^2}{S}$,
we have $\widetilde H\ge\mathcal E_0$.
We collect the remaining surface tension and vorticity terms in
\begin{align}
R_\kappa={}&-\frac\kappa2\int_{-\pi}^{\pi}\frac{m_{xx}\eta\eta_x}{S}\d x+\kappa\int_{-\pi}^{\pi}m_x\left(1-\frac{1+\eta_x^2/2}{S}\right)\d x\notag\\
&+\frac{3\kappa}{4}\int_{-\pi}^{\pi}\frac{\eta_x^2}{S}\d x+g\int_{-\pi}^{\pi}\chi\eta^2\d x,\label{Rkappa}\\
R_\gamma={}&\gamma\int_{-\pi}^{\pi}\zeta\partial_x^{-1}w\d x+\gamma\int_{-\pi}^{\pi}\left[\frac12m_{xx}\eta^2-\left(\frac54+\frac12m_x\right)\eta\eta_x\right]\psi\d x.
\label{Rgamma}
\end{align}

Combining \eqref{potential-identity} and \eqref{equipartition} gives the following identity.

\begin{lemma}\label{Lemma-combined}
Every regular solution of \eqref{Waterwave} satisfies
\begin{equation}\label{combined-identity}
\begin{aligned}
\frac12\widetilde H+R_\kappa={}&\frac32\int_{-\pi}^{\pi}\chi\psi w\d x-\frac{\d}{\d t}\int_{-\pi}^{\pi}\zeta\psi\d x-\int_{-\pi}^{\pi}P_{\rm ext}\zeta\d x\\
&-\int_{-\pi}^{\pi}m\psi_xw\d x-\int_{-\pi}^{\pi}\zeta N(\eta)\psi\d x+R_\gamma.
\end{aligned}
\end{equation}
\end{lemma}

\begin{proof}
In \eqref{equipartition}, take
\[
\vartheta=\frac32\chi-\frac14=\frac54-\frac32m_x.
\]
Add this identity to \eqref{potential-identity} after moving its $\chi$-weighted potential term to the left. The terms involving time derivatives and pressure combine because $\partial_x(m\eta)+\vartheta\eta=\zeta$.  On the right,
\[
\int_{-\pi}^{\pi}\vartheta\psi w\d x=\frac32\int_{-\pi}^{\pi}\chi\psi w\d x-\frac14\int_{-\pi}^{\pi}\psi w\d x.
\]
Move the last term to the left. The coefficient of $g\eta^2$ there is $m_x/2+\vartheta=5/4-m_x$. The coefficient of $\kappa\eta_x^2/S$, apart from $R_1$, is $3m_x/2+\vartheta=5/4$. Since $\vartheta_x=-3m_{xx}/2$, the mixed capillary term is $-\kappa\int_{-\pi}^{\pi}m_{xx}\eta\eta_x/S\d x/2$. These are exactly the terms in $\widetilde H/2+R_\kappa$. The nonlocal vorticity terms combine into $\gamma\int_{-\pi}^{\pi}\zeta\partial_x^{-1}w\d x$. Put $b=5/4-m_x/2$ for this calculation. The local vorticity terms are
\[
\gamma\int_{-\pi}^{\pi}\bigl(b\eta^2\psi_x+\vartheta\psi\eta\eta_x\bigr)\d x.
\]
Integration by parts, using $b_x=-m_{xx}/2$ and $\vartheta-2b=-5/4-m_x/2$, gives the second term in \eqref{Rgamma}. This proves \eqref{combined-identity}.
\end{proof}

\subsection{Periodic Pohozaev identities and the energy inequality}
The last two boundary integrals in \eqref{combined-identity} have no definite sign. We use the harmonicity of $\Phi$ to express them as interior integrals and terms on the bed. We shall use the periodic form of the Rellich identity appearing in \cite[Appendix~B]{Alazard17}, together with two Pohozaev identities. For the bed integral, write
\[
B_0(t)=\int_{-\pi}^{\pi}\Phi_x(t,x,-h)^2\d x.
\]

\begin{lemma}\label{Lemma-rellich}
For every smooth periodic function $\mu=\mu(x)$,
\begin{equation}\label{rellich}
\int_{-\pi}^{\pi}\mu N(\eta)\psi\d x=\frac12\int_{-\pi}^{\pi}\mu\Phi_x(x,-h)^2\d x-\int_{D_\eta}\mu_x\Phi_x\Phi_y\d x\d y.
\end{equation}
Also,
\begin{align}
\int_{-\pi}^{\pi}m\psi_xw\d x+\int_{-\pi}^{\pi}m\eta_xN(\eta)\psi\d x&=\frac12\int_{D_\eta}m_x(\Phi_x^2-\Phi_y^2)\d x\d y,\label{periodic-pohozaev}\\
\frac12\int_{D_\eta}(\Phi_x^2-\Phi_y^2)\d x\d y&=\int_{-\pi}^{\pi}\eta N(\eta)\psi\d x+\frac h2B_0.\label{vertical-pohozaev}
\end{align}
\end{lemma}

\begin{proof}
We first consider
\[
X=\left(-\mu\Phi_x\Phi_y,\frac{\mu}{2}(\Phi_x^2-\Phi_y^2)\right).
\]
The harmonicity of $\Phi$ gives $\operatorname{div}X=-\mu_x\Phi_x\Phi_y$. Integrating over $D_\eta$ and using \eqref{Nbulk} and the bed condition yields \eqref{rellich}.

For \eqref{periodic-pohozaev}, we use
\[
Y=\left(\frac{m}{2}(\Phi_x^2-\Phi_y^2),m\Phi_x\Phi_y\right),
\]
whose divergence is $m_x(\Phi_x^2-\Phi_y^2)/2$. On the free surface, \eqref{Nbulk} gives
\[
Y\cdot(-\eta_x,1)=m\bigl(\psi_xw+\eta_xN(\eta)\psi\bigr).
\]
Since the bed integral vanishes, the divergence theorem gives \eqref{periodic-pohozaev}.
Finally, taking
\[
Z=\left(-y\Phi_x\Phi_y,\frac{y}{2}(\Phi_x^2-\Phi_y^2)\right),
\]
we have $\operatorname{div}Z=(\Phi_x^2-\Phi_y^2)/2$. The free-surface and bed integrals are $\int_{-\pi}^{\pi}\eta N(\eta)\psi\d x$ and $hB_0/2$, respectively. Thus the same argument proves \eqref{vertical-pohozaev}.
\end{proof}

To combine the boundary terms in \eqref{combined-identity}, set
\begin{equation}
\rho=\zeta+\eta-m\eta_x=\left(\frac94-\frac12m_x\right)\eta.
\end{equation}
The preceding identities give the following form of \eqref{combined-identity}.

\begin{lemma}\label{Lemma-pohozaev-combined}
Every regular solution of \eqref{Waterwave} satisfies
\begin{equation}\label{pohozaev-combined}
\begin{aligned}
\frac12\widetilde H+R_\kappa+\frac12\int_{-\pi}^{\pi}(h+\rho)\Phi_x(x,-h)^2\d x &=\frac32\int_{-\pi}^{\pi}\chi\psi w\d x+\frac12\int_{D_\eta}\chi(\Phi_x^2-\Phi_y^2)\d x\d y\\
&\quad-\frac{\d}{\d t}\int_{-\pi}^{\pi}\zeta\psi\d x-\int_{-\pi}^{\pi}P_{\rm ext}\zeta\d x\\
&\quad+R_\gamma+\int_{D_\eta}\rho_x\Phi_x\Phi_y\d x\d y.
\end{aligned}
\end{equation}
\end{lemma}

\begin{proof}
We first use \eqref{periodic-pohozaev} to combine the two terms containing $m$:
\[
-\int_{-\pi}^{\pi}m\psi_xw\d x-\int_{-\pi}^{\pi}\zeta N(\eta)\psi\d x=-\frac12\int_{D_\eta}m_x(\Phi_x^2-\Phi_y^2)\d x\d y-\int_{-\pi}^{\pi}(\zeta-m\eta_x)N(\eta)\psi\d x.
\]
Since $m_x=1-\chi$, the first integral on the right can be split into an unweighted integral and one containing $\chi$. Applying \eqref{vertical-pohozaev} to the unweighted integral, and using $\rho=\zeta+\eta-m\eta_x$, we obtain
\[
-\int_{-\pi}^{\pi}m\psi_xw\d x-\int_{-\pi}^{\pi}\zeta N(\eta)\psi\d x=\frac12\int_{D_\eta}\chi(\Phi_x^2-\Phi_y^2)\d x\d y-\frac h2 B_0-\int_{-\pi}^{\pi}\rho N(\eta)\psi\d x.
\]
It remains to express the last boundary integral in terms of $\Phi$. Taking $\mu=\rho$ in \eqref{rellich} gives
\[
\int_{-\pi}^{\pi}\rho N(\eta)\psi\d x=\frac12\int_{-\pi}^{\pi}\rho\Phi_x(x,-h)^2\d x-\int_{D_\eta}\rho_x\Phi_x\Phi_y\d x\d y.
\]
Combining these identities gives
\[
\begin{aligned}
-\int_{-\pi}^{\pi}m\psi_xw\d x-\int_{-\pi}^{\pi}\zeta N(\eta)\psi\d x&=\frac12\int_{D_\eta}\chi(\Phi_x^2-\Phi_y^2)\d x\d y-\frac12\int_{-\pi}^{\pi}(h+\rho)\Phi_x(x,-h)^2\d x\\
&\quad+\int_{D_\eta}\rho_x\Phi_x\Phi_y\d x\d y.
\end{aligned}
\]
Substituting this expression into \eqref{combined-identity} and moving the bed integral to the left yields
\eqref{pohozaev-combined}.
\end{proof}

We now express \eqref{pohozaev-combined} in terms of $\widetilde\psi$, using the surface equation to account for the spatial mean of $\psi$. Set
\[
\nu(t)=\int_{-\pi}^{\pi}\chi\eta\d x,\quad \beta(t)=\frac{3}{8\pi}\nu(t),\quad I(t)=\int_{-\pi}^{\pi}\zeta\widetilde\psi\d x.
\]
Since $\langle\eta\rangle=0$,
\begin{equation}\label{zeta-mean}
\int_{-\pi}^{\pi}\zeta\d x=\frac32\nu(t).
\end{equation}
After separating the mean of $\psi$, the vorticity terms take the form
\begin{equation}\label{Rgamma-mean}
\widetilde R_\gamma=\gamma\int_{-\pi}^{\pi}\zeta\partial_x^{-1}w\d x+\gamma\int_{-\pi}^{\pi}\left[\frac12m_{xx}\eta^2-\left(\frac54+\frac12m_x\right)\eta\eta_x\right]\widetilde\psi\d x-\gamma\langle\eta\widetilde\psi_x\rangle\int_{-\pi}^{\pi}\zeta\d x.
\end{equation}
With this notation, \eqref{pohozaev-combined} takes the following form.

\begin{prop}\label{Prop-multiplier}
Every regular solution of \eqref{Waterwave} satisfies
\begin{equation}\label{multiplier-identity}
\begin{aligned}
&\frac12\widetilde H+R_\kappa+\int_{-\pi}^{\pi}\left(\frac{h+\rho}{2}-\beta\right)\Phi_x(x,-h)^2\d x\\
&\quad=\frac32\int_{-\pi}^{\pi}\chi\widetilde\psi w\d x+\frac12\int_{D_\eta}\chi(\Phi_x^2-\Phi_y^2)\d x\d y-I'(t)\\
&\qquad-\int_{-\pi}^{\pi}\widetilde P_{\rm ext}\zeta\d x+\widetilde R_\gamma+\int_{D_\eta}\rho_x\Phi_x\Phi_y\d x\d y.
\end{aligned}
\end{equation}
Here $\widetilde P_{\rm ext}=P_{\rm ext}-\langle P_{\rm ext}\rangle$.
\end{prop}

\begin{proof}
Put $p(t)=\langle\psi\rangle$. Taking $\mu=1$ in \eqref{rellich} gives
\begin{equation}\label{mean-N}
\langle N(\eta)\psi\rangle=\frac{B_0}{4\pi}.
\end{equation}
Taking the mean of the second equation in \eqref{Waterwave}, we obtain
\[
p'=-\frac{B_0}{4\pi}+\gamma\langle\eta\widetilde\psi_x\rangle-\langle P_{\rm ext}\rangle.
\]
In \eqref{pohozaev-combined}, replace $\psi$ by $\widetilde\psi+p$. The change in the local part of $R_\gamma$ is
\[
R_\gamma(\eta,\psi)-R_\gamma(\eta,\widetilde\psi)=\frac{3\gamma p}{4}\int_{-\pi}^{\pi}m_{xx}\eta^2\d x.
\]
On the other hand,
\begin{align*}
\frac32p\int_{-\pi}^{\pi}\chi w\d x-p\frac{\d}{\d t}\int_{-\pi}^{\pi}\zeta\d x=-\frac{3\gamma p}{2}\int_{-\pi}^{\pi}\chi\eta\eta_x\d x=-\frac{3\gamma p}{4}\int_{-\pi}^{\pi}m_{xx}\eta^2\d x.
\end{align*}
Thus these terms cancel. The remaining mean contribution is
\[
-p'\int_{-\pi}^{\pi}\zeta\d x={}\beta B_0-\gamma\langle\eta\widetilde\psi_x\rangle\int_{-\pi}^{\pi}\zeta\d x+\langle P_{\rm ext}\rangle\int_{-\pi}^{\pi}\zeta\d x.
\]
The first term is moved to the left, the second is included in \eqref{Rgamma-mean}, and the last combines with the pressure integral. This proves \eqref{multiplier-identity}.
\end{proof}

We finally estimate the geometric terms in \eqref{multiplier-identity}. The condition \eqref{capillary-cutoff-condition} allows the mixed surface tension term to be controlled by the gravitational and surface energies. 
For $\|\eta\|_{W^{2,\infty}}$ sufficiently small, the coefficient of the bed integral is positive, and the interior term containing $\rho_x$ is absorbed into the kinetic energy.

\begin{lemma}\label{Lemma-geometry}
There is $\delta_2>0$, depending on $g,\kappa,h$ and $m$, such that, if
$\|\eta\|_{W^{2,\infty}}\le\delta_2$, then
\begin{equation}\label{geometry-bound}
\begin{aligned}
\frac12\widetilde H+R_\kappa+\int_{-\pi}^{\pi}\left(\frac{h+\rho}{2}-\beta\right)\Phi_x(x,-h)^2\d x-\int_{D_\eta}\rho_x\Phi_x\Phi_y\d x\d y\ge\frac18\mathcal E_0+\frac h4B_0.
\end{aligned}
\end{equation}
\end{lemma}

\begin{proof}
We may take $\|\eta_x\|_{L^\infty}\le1$. Young's inequality and \eqref{capillary-cutoff-condition} give
\begin{align*}
\left|\frac\kappa2\int_{-\pi}^{\pi}
\frac{m_{xx}\eta\eta_x}{S}\d x\right|&\le\frac\kappa8\int_{-\pi}^{\pi}\frac{m_{xx}^2\eta^2}{S}\d x+\frac\kappa2\int_{-\pi}^{\pi}\frac{\eta_x^2}{S}\d x\\
&\le\frac g8\int_{-\pi}^{\pi}\eta^2\d x+\frac\kappa2\int_{-\pi}^{\pi}\frac{\eta_x^2}{S}\d x.
\end{align*}
Moreover,
\[
1-\frac{1+\eta_x^2/2}{S}=-\frac{(S-1)^2}{2S},\quad\frac{(S-1)^2}{2S}\le\frac{\eta_x^4}{8S}.
\]
Decrease $\delta_2$ so that $\|m_x\|_{L^\infty}\|\eta_x\|_{L^\infty}^2\le1$. The second term in \eqref{Rkappa} is then bounded in absolute value by $\kappa\int_{-\pi}^{\pi}\eta_x^2/S\d x/8$. Since $g\int_{-\pi}^{\pi}\chi\eta^2\d x\ge0$, we obtain
\[
\frac12\widetilde H+R_\kappa\ge\frac12K+\frac g8\int_{-\pi}^{\pi}\eta^2\d x+\frac{5\kappa}{8}\int_{-\pi}^{\pi}\frac{\eta_x^2}{S}\d x\ge\frac14\mathcal E_0.
\]
Also $\|\rho\|_{W^{1,\infty}}+|\beta|\le C\delta_2$. We can therefore arrange
\[
\frac{h+\rho}{2}-\beta\ge\frac h4,\quad\left|\int_{D_\eta}\rho_x\Phi_x\Phi_y\d x\d y\right|\le\frac18K.
\]
These estimates prove \eqref{geometry-bound}.
\end{proof}

The preceding lemma provides the lower bound needed to turn the multiplier identity into an estimate for the time integral of the energy. Let $0\le a<b\le T$ be such that the smallness condition of Lemma~\ref{Lemma-geometry} holds throughout $[a,b]$. Integrating \eqref{multiplier-identity} over this interval and using \eqref{geometry-bound}, we obtain
\begin{equation}\label{fundamental-inequality}
\frac18\int_a^b\mathcal E_0(t)\d t+\frac h4\int_a^b B_0(t)\d t\le\mathcal O+\mathcal W+\mathcal B+\int_a^b\widetilde R_\gamma(t)\d t,
\end{equation}
where
\begin{align}
\mathcal O={}&\int_a^b\left[\frac32\int_{-\pi}^{\pi}\chi\widetilde\psi w\d x+\frac12\int_{D_\eta}\chi(\Phi_x^2-\Phi_y^2)\d x\d y\right]\d t,\label{Odef}\\
\mathcal W={}&-\int_a^b\int_{-\pi}^{\pi}\widetilde P_{\rm ext}\zeta\d x\d t,\quad\mathcal B=-[I(t)]_a^b.\label{WBdef}
\end{align}
Here $\mathcal O$ contains a surface integral over $x\in\omega_{\rm abs}$ and an interior integral over the fluid below that part of the free surface, both weighted by $\chi$. The term $\mathcal W$ is the pressure contribution, and $\mathcal B$ consists of the values of $I$ at the time endpoints.

\section{Exponential decay of the energy}\label{Section-decay}

\subsection{Estimates of the pressure, vorticity and observation terms}

We use the feedback law \eqref{Feedback-law} to estimate the right-hand side of \eqref{fundamental-inequality} on an interval $[a,b]\subset[0,T]$ where the hypotheses of Lemmas~\ref{Lemma-elliptic}, \ref{Lemma-coercivity} and \ref{Lemma-geometry} hold. Throughout this subsection, constants depend on $g,\kappa,h,\gamma$ and the fixed cut-off, with any additional dependence indicated by subscripts, and are independent of the solution and of $a,b$.

We begin by bounding the pressure in terms of the energy loss and expressing the contribution of its mean as an endpoint term.

\begin{lemma}\label{Lemma-pressure}
Every regular solution of \eqref{Waterwave} with feedback \eqref{Feedback-law} satisfies
\begin{equation}\label{pressure-L2}
\int_a^b\int_{-\pi}^{\pi}P_{\rm ext}^2\d x\d t\le\|\chi\|_{L^\infty}\bigl(\mathcal E_\gamma(a)-\mathcal E_\gamma(b)\bigr).
\end{equation}
Moreover,
\begin{equation}\label{pressure-mean-endpoint}
\int_a^b\langle P_{\rm ext}\rangle\left(\int_{-\pi}^{\pi}\zeta\d x\right)\d t=\frac{3}{8\pi}[\nu(t)^2]_a^b.
\end{equation}
\end{lemma}

\begin{proof}
Since $P_{\rm ext}=\chi\eta_t$ and $\chi\ge0$, we have
\[
\int_{-\pi}^{\pi}P_{\rm ext}^2\d x=\int_{-\pi}^{\pi}\chi^2\eta_t^2\d x\le\|\chi\|_{L^\infty}\mathcal D(t).
\]
Integration in time, followed by the energy identity \eqref{dissipation}, gives \eqref{pressure-L2}.

The cut-off is independent of time. Hence
\[
\nu'(t)=\int_{-\pi}^{\pi}\chi\eta_t\d x=\int_{-\pi}^{\pi}P_{\rm ext}\d x,\quad\langle P_{\rm ext}\rangle=\frac{\nu'(t)}{2\pi}.
\]
Using \eqref{zeta-mean}, we obtain
\[
\int_a^b\langle P_{\rm ext}\rangle\left(\int_{-\pi}^{\pi}\zeta\d x\right)\d t=\frac{3}{4\pi}\int_a^b\nu'(t)\nu(t)\d t=\frac{3}{8\pi}[\nu(t)^2]_a^b.
\]
This proves \eqref{pressure-mean-endpoint}.
\end{proof}

The pressure mean thus contributes an endpoint term, which can be estimated together with $\mathcal B$ by the energy at $t=a,b$. Using the feedback law to estimate the remaining pressure integral, we obtain the following bound for $\mathcal W+\mathcal B$ in \eqref{fundamental-inequality}.

\begin{lemma}\label{Lemma-pressure-work}
For every $\epsilon>0$, there is $C_\epsilon>0$ such that
\begin{equation}\label{pressure-work-bound}
\mathcal W+\mathcal B\le C_\epsilon\bigl(\mathcal E_\gamma(a)+\mathcal E_\gamma(b)\bigr) +C_\epsilon\int_a^b\mathcal D(t)\d t +\epsilon\int_a^b\mathcal E_0(t)\d t.
\end{equation}
\end{lemma}

\begin{proof}
The definition of $\zeta$ and the boundedness of $m,m_x$ imply
\[
\|\zeta\|_{L^2}\le C\bigl(\|\eta_x\|_{L^2}+\|\eta\|_{L^2}\bigr)\le C\|\eta\|_{H^1}.
\]
Consequently, by \eqref{Knorm} and the energy equivalence for
$\mathcal E_0$,
\[
|I(t)|\le\|\zeta\|_{L^2}\|\widetilde\psi\|_{L^2}\le C\mathcal E_0(t),\quad|\nu(t)|^2\le\|\chi\|_{L^2}^2\|\eta(t)\|_{L^2}^2\le C\mathcal E_0(t).
\]
The coercivity estimate therefore controls these quantities by
$C\mathcal E_\gamma(t)$.

Separating the mean of the pressure in \eqref{WBdef} and using \eqref{pressure-mean-endpoint}, we find
\begin{equation}\label{work-split}
\mathcal W=-\int_a^b\int_{-\pi}^{\pi}P_{\rm ext}\zeta\d x\d t+\frac{3}{8\pi}[\nu^2]_a^b.
\end{equation}
The endpoint term in \eqref{work-split} and $\mathcal B$ are bounded by $C(\mathcal E_\gamma(a)+\mathcal E_\gamma(b))$. For the remaining integral, the Cauchy--Schwarz inequality gives
\[
\begin{aligned}
\left|\int_{-\pi}^{\pi}P_{\rm ext}\zeta\d x\right|&\le\left(\int_{-\pi}^{\pi}\chi\eta_t^2\d x\right)^{1/2}\left(\int_{-\pi}^{\pi}\chi\zeta^2\d x\right)^{1/2}\le C\mathcal D(t)^{1/2}\mathcal E_0(t)^{1/2}\\
&\le\epsilon\mathcal E_0(t)+C_\epsilon\mathcal D(t).
\end{aligned}
\]
Integrating this estimate over $[a,b]$ and combining it with the preceding bounds for the endpoint terms, we obtain \eqref{pressure-work-bound}.
\end{proof}

We next estimate $\widetilde R_\gamma$ in \eqref{Rgamma-mean}. Using the kinematic condition to extract a time derivative, we obtain the following bound.

\begin{lemma}\label{Lemma-gamma-endpoint}
For every $\epsilon>0$, there are $\delta_\epsilon>0$ and $C_{\epsilon,\gamma}>0$ such that, if
\[
\sup_{a\le t\le b}\|\eta(t)\|_{W^{2,\infty}} \le\delta_\epsilon,
\]
then
\begin{equation}\label{gamma-integral-bound}
\int_a^b\widetilde R_\gamma(t)\d t\le{}C_\gamma\bigl(\mathcal E_\gamma(a)+\mathcal E_\gamma(b)\bigr)+\epsilon\int_a^b\mathcal E_0(t)\d t+C_{\epsilon,\gamma}\int_a^b\|\widetilde\psi(t)\|_{L^2}^2\d t.
\end{equation}
Here $C_\gamma$ is independent of $\epsilon$.
\end{lemma}

\begin{proof}
With $\vartheta$ chosen as in the proof of Lemma~\ref{Lemma-combined}, we integrate by parts in the term containing $\partial_x(m\eta)$ and use the kinematic condition to obtain
\[
\begin{aligned}
\gamma\int_{-\pi}^{\pi}\zeta\partial_x^{-1}w\d x&=\gamma\int_{-\pi}^{\pi}\vartheta\eta\partial_x^{-1}w\d x-\gamma\int_{-\pi}^{\pi}m\eta w\d x\\
&=-\frac{\gamma}{2}\frac{\d}{\d t}\int_{-\pi}^{\pi}m\eta^2\d x+\gamma\int_{-\pi}^{\pi}\vartheta\eta\partial_x^{-1}w\d x-\frac{\gamma^2}{3}\int_{-\pi}^{\pi}m_x\eta^3\d x.
\end{aligned}
\]
Substituting this expression into \eqref{Rgamma-mean} gives
\begin{equation}\label{gamma-endpoint}
\widetilde R_\gamma=-\frac{\gamma}{2}\frac{\d}{\d t}\int_{-\pi}^{\pi}m\eta^2\d x+\gamma\int_{-\pi}^{\pi}\vartheta\eta\partial_x^{-1}w\d x+Q_\gamma,
\end{equation}
where
\begin{equation}\label{Qgamma}
\begin{aligned}
Q_\gamma={}&-\frac{\gamma^2}{3}\int_{-\pi}^{\pi}m_x\eta^3\d x+\gamma\int_{-\pi}^{\pi}\left[\frac12m_{xx}\eta^2-\left(\frac54+\frac12m_x\right)\eta\eta_x\right]\widetilde\psi\d x\\
&-\gamma\langle\eta\widetilde\psi_x\rangle\int_{-\pi}^{\pi}\zeta\d x.
\end{aligned}
\end{equation}
We first estimate $Q_\gamma$ in terms of $\mathcal E_0$. By the Cauchy--Schwarz inequality and the boundedness of $m_x,m_{xx}$, the absolute values of the first two integrals in \eqref{Qgamma} are bounded by
\[
C\gamma^2\|\eta\|_{L^\infty}\|\eta\|_{L^2}^2+C|\gamma|\|\eta\|_{L^\infty}\|\eta\|_{H^1}\|\widetilde\psi\|_{L^2}.
\]
For the last term, periodic integration by parts transfers the derivative from $\widetilde\psi$ to $\eta$. Together with $\int_{-\pi}^{\pi}\zeta\d x=\frac32\int_{-\pi}^{\pi}\chi\eta\d x$, this gives
\[
\left|\gamma\langle\eta\widetilde\psi_x\rangle\int_{-\pi}^{\pi}\zeta\d x\right|=|\gamma|\left|\langle\eta_x\widetilde\psi\rangle\right|\left|\int_{-\pi}^{\pi}\zeta\d x\right|\leq C|\gamma|\|\eta\|_{L^\infty}
\|\eta_x\|_{L^2}\|\widetilde\psi\|_{L^2}.
\]
The energy equivalence for $\mathcal E_0$ therefore yields
\begin{equation}\label{Qgamma-bound}
|Q_\gamma|\le C(|\gamma|+\gamma^2)\|\eta\|_{W^{2,\infty}}\mathcal E_0.
\end{equation}

It remains to control the first two terms in \eqref{gamma-endpoint}. After integration in time, the derivative contributes endpoint values, which are bounded by the energy since \eqref{coercivity} gives
\[
\left|\frac{\gamma}{2}\int_{-\pi}^{\pi}m\eta^2\d x\right|\le C_\gamma\|\eta\|_{L^2}^2\le C_\gamma\mathcal E_\gamma(t).
\]
For the nonlocal integral, we use \eqref{AGbound} followed by Young's inequality to obtain, for every $\epsilon_1>0$,
\begin{equation}\label{gamma-low-order}
\left|\gamma\int_{-\pi}^{\pi}\vartheta\eta\partial_x^{-1}w\d x\right|\le C|\gamma|\|\eta\|_{L^2}\|\widetilde\psi\|_{L^2}\le\epsilon_1\mathcal E_0+C_{\epsilon_1,\gamma}\|\widetilde\psi\|_{L^2}^2.
\end{equation}
Taking $\epsilon_1=\epsilon/2$ and decreasing $\delta_\epsilon$ so that $C(|\gamma|+\gamma^2)\delta_\epsilon\le\epsilon/2$,
we see from \eqref{Qgamma-bound} and \eqref{gamma-low-order} that the last two terms in \eqref{gamma-endpoint} are bounded above by $\epsilon\mathcal E_0+C_{\epsilon,\gamma}\|\widetilde\psi\|_{L^2}^2$.
Integrating \eqref{gamma-endpoint} over $[a,b]$ and combining this estimate with the endpoint bound proves \eqref{gamma-integral-bound}.
\end{proof}

We next estimate $\mathcal O$ in \eqref{Odef}, using Green's formula to relate its interior integral to the surface integral.

\begin{lemma}\label{Lemma-observation_terms}
For every $\epsilon>0$, there are $\delta_\epsilon>0$ and $C_\epsilon>0$ such that, if
\[
\sup_{a\le t\le b}\|\eta(t)\|_{W^{2,\infty}}\le\delta_\epsilon,
\]
then
\begin{equation}\label{observation-bound}
|\mathcal O|\le\epsilon\int_a^b\mathcal E_0(t)\d t+C_\epsilon\int_a^b\mathcal D(t)\d t+C_\epsilon\int_a^b\|\widetilde\psi(t)\|_{L^2}^2\d t.
\end{equation}
\end{lemma}

\begin{proof}
Using the kinematic condition, we write
\[
\int_{-\pi}^{\pi}\chi\widetilde\psi w\d x=\int_{-\pi}^{\pi}\chi\widetilde\psi\eta_t\d x-\gamma\int_{-\pi}^{\pi}\chi\widetilde\psi\eta\eta_x\d x.
\]
By \eqref{Knorm} and the Cauchy--Schwarz and Young inequalities, for every $\epsilon_1>0$ we have
\begin{equation}\label{surface-observation-bound}
\begin{aligned}
\left|\int_{-\pi}^{\pi}\chi\widetilde\psi w\d x\right|&\le C\mathcal E_0^{1/2}\mathcal D^{1/2}+C|\gamma|\delta_\epsilon\|\widetilde\psi\|_{L^2}\|\eta_x\|_{L^2}\\
&\le(\epsilon_1+C|\gamma|\delta_\epsilon)\mathcal E_0+C_{\epsilon_1}\mathcal D.
\end{aligned}
\end{equation}
For the interior integral, applying the divergence theorem to $\chi\widetilde\Phi\nabla\Phi$ and using the boundary conditions gives
\begin{equation*}
\int_{D_\eta}\chi|\nabla\Phi|^2\d x\d y=\int_{-\pi}^{\pi}\chi\widetilde\psi w\d x-\int_{D_\eta}\chi_x\widetilde\Phi\Phi_x\d x\d y.
\end{equation*}
Since $\chi\ge0$, we can use \eqref{PhiL2} to obtain, for every
$\epsilon_2>0$,
\[
\begin{aligned}
\left|\int_{D_\eta}\chi(\Phi_x^2-\Phi_y^2)\d x\d y\right|&\le\int_{D_\eta}\chi|\nabla\Phi|^2\d x\d y\\
&\le\left|\int_{-\pi}^{\pi}\chi\widetilde\psi w\d x\right|+C\|\widetilde\Phi\|_{L^2(D_\eta)}K^{1/2}\\
&\le\left|\int_{-\pi}^{\pi}\chi\widetilde\psi w\d x\right|+\epsilon_2K+C_{\epsilon_2}\|\widetilde\psi\|_{L^2}^2.
\end{aligned}
\]
Combining this with \eqref{surface-observation-bound}, we bound the integrand in \eqref{Odef} in absolute value by
\[
\left(2\epsilon_1+2C|\gamma|\delta_\epsilon+\frac{\epsilon_2}{2}\right)\mathcal E_0+2C_{\epsilon_1}\mathcal D+\frac{C_{\epsilon_2}}{2}\|\widetilde\psi\|_{L^2}^2.
\]
Choosing $2\epsilon_1+\epsilon_2/2\le\epsilon/2$ and then decreasing $\delta_\epsilon$ so that $2C|\gamma|\delta_\epsilon\le\epsilon/2$, we obtain \eqref{observation-bound} after integration over $[a,b]$.
\end{proof}

Combining the preceding estimates with \eqref{fundamental-inequality} gives the following energy bound with an $L^2$ remainder.

\begin{prop}\label{Prop-weak}
There are $\delta_*>0$ and $C_w\ge1$ such that every regular solution of \eqref{Waterwave} with feedback \eqref{Feedback-law} for which
\[
\sup_{a\le t\le b}\|\eta(t)\|_{W^{2,\infty}}\le\delta_*
\]
satisfies the estimate
\begin{equation}\label{weak-estimate}
\int_a^b\mathcal E_\gamma(t)\d t\le C_w\Bigg[\mathcal E_\gamma(a)+\mathcal E_\gamma(b)+\int_a^b\mathcal D(t)\d t+\int_a^b\|\widetilde\psi(t)\|_{L^2}^2\d t\Bigg].
\end{equation}
The constant $C_w$ is independent of $a,b$.
\end{prop}

\begin{proof}
We take $\epsilon=1/64$ in \eqref{pressure-work-bound}, \eqref{gamma-integral-bound} and \eqref{observation-bound}, and choose $\delta_*$ so that these estimates, \eqref{geometry-bound} and \eqref{coercivity} apply. Substituting into \eqref{fundamental-inequality} and moving the energy remainder to the left gives
\[
\frac{5}{64}\int_a^b\mathcal E_0(t)\d t+\frac h4\int_a^b B_0(t)\d t\le{}C\Bigg[\mathcal E_\gamma(a)+\mathcal E_\gamma(b)+\int_a^b\mathcal D(t)\d t+\int_a^b\|\widetilde\psi(t)\|_{L^2}^2\d t\Bigg].
\]
Dropping the nonnegative bed integral and using $\mathcal E_\gamma\le3\mathcal E_0/2$, we obtain \eqref{weak-estimate} with $C_w\ge1$.
\end{proof}

\subsection{Unique continuation for the linearized equations}

To remove the $L^2$ remainder in \eqref{weak-estimate}, we use a compactness argument based on unique continuation (cf. \cite{DLZ03}). We first establish the required unique continuation property for the linearized equations.

Set $G_0=G(0)$, whose symbol on nonzero Fourier modes is $|j|\tanh(h|j|)$. Thus, for every $r\in\mathbb R$, the operators
\[
G_0:H^{r+1}_0\rightarrow H^r_0,\quad g-\kappa\partial_x^2:H^{r+2}\rightarrow H^r
\]
are isomorphisms. Linearizing \eqref{Waterwave} about the flat Couette flow with $P_{\rm ext}=0$, we obtain
\begin{equation}\label{linear-system}
\begin{cases}
\partial_t\eta_{\rm lin}=G_0\psi_{\rm lin},\\ \partial_t\psi_{\rm lin}=-(g-\kappa\partial_x^2)\eta_{\rm lin}+\gamma\partial_x^{-1}G_0\psi_{\rm lin},
\end{cases}
\end{equation}
where $\eta_{\rm lin}$ and $\psi_{\rm lin}$ are the zero-mean linearized surface variables. We shall use the following consequence of harmonicity.

\begin{lemma}\label{Lemma-harmonic-unique}
Let $q$ be a smooth periodic function, and let $Q$ solve
\[
\Delta Q=0\ (-h<y<0),\quad Q|_{y=0}=q,\quad Q_y|_{y=-h}=0.
\]
If
\[
q_x=0,\quad G_0q=0
\]
on a nonempty open interval $J$ of the upper boundary, then $q$ is constant. In particular, $q=0$ if it has zero mean.
\end{lemma}

\begin{proof}
Since $q_x=0$ on $J$, we have $q=c$ there for some constant $c$, while $G_0q=0$ gives $Q_y=0$ on the same segment. We extend $Q-c$ by zero across a smaller segment of $J$. Both boundary contributions vanish in Green's formula, so the extension is harmonic in the sense of distributions and hence is smooth and harmonic. Since it vanishes on one side of the boundary, unique continuation gives $Q=c$ near the segment and then throughout the strip. Thus $q=c$, with $c=0$ when $q$ has zero mean.
\end{proof}

We now apply this lemma to time derivatives of the linearized solution.

\begin{prop}\label{Prop-linear-unique}
Let $(\eta_{\rm lin},\psi_{\rm lin})$ solve \eqref{linear-system} in distributions on $(0,T)$, where $T>0$, and suppose that
\[
(\eta_{\rm lin},\psi_{\rm lin})\in L^\infty(0,T;H^1_0\times H^{1/2}_0).
\]
If $\partial_t\eta_{\rm lin}=0$ in distributions on $(0,T)\times J$ for a nonempty open interval $J$, then $\eta_{\rm lin}=\psi_{\rm lin}=0$.
\end{prop}

\begin{proof}
We first consider a smooth solution. Differentiating \eqref{linear-system} and using $\partial_t\eta_{\rm lin}=0$ on $(0,T)\times J$, we obtain
\[
\begin{aligned}
G_0\partial_t^2\psi_{\rm lin}&=\partial_t^3\eta_{\rm lin}=0,\\
\partial_x\partial_t^2\psi_{\rm lin}&=-(g-\kappa\partial_x^2)\partial_x\partial_t\eta_{\rm lin}+\gamma \partial_t^2\eta_{\rm lin}=0
\end{aligned}\qquad\text{on }(0,T)\times J.
\]
Applying Lemma~\ref{Lemma-harmonic-unique} to $\partial_t^2\psi_{\rm lin}$ at each time and using its zero mean, we conclude that $\partial_t^2\psi_{\rm lin}=0$ throughout the period.
Differentiating the second equation in \eqref{linear-system} twice and once in time, respectively, now gives
\[
(g-\kappa\partial_x^2)\partial_t^2\eta_{\rm lin}=0,\quad(g-\kappa\partial_x^2)\partial_t\eta_{\rm lin}=\gamma\partial_x^{-1}\partial_t^2\eta_{\rm lin}.
\]
The invertibility of $g-\kappa\partial_x^2$ yields $\partial_t^2\eta_{\rm lin}=0$ and then $\partial_t\eta_{\rm lin}=0$. Hence $G_0\psi_{\rm lin}=0$, and the zero mean of $\psi_{\rm lin}$ implies $\psi_{\rm lin}=0$. The second equation then gives $\eta_{\rm lin}=0$.

For regular solution, one can convolve in time with a smooth mollifier of radius $0<\varepsilon<T/2$ and denote the result by $(\eta_{\rm lin}^{\varepsilon},\psi_{\rm lin}^{\varepsilon})$ on $(\varepsilon,T-\varepsilon)$. Since the coefficients are independent of time, this pair satisfies \eqref{linear-system} and $\partial_t\eta_{\rm lin}^{\varepsilon}=0$ on $J$, and is smooth in time with values in $H^1_0\times H^{1/2}_0$.
The equations also improve its spatial regularity. More precisely, it follows that for every integer
$k\ge 0$,
\[
\partial_t^k \psi_{\rm lin}^{\varepsilon}=G_0^{-1}\partial_t^{k+1}\eta_{\rm lin}^{\varepsilon}\in H^2_0,
\]
and consequently,
\[
(g-\kappa\partial_x^2)\partial_t^k \eta_{\rm lin}^{\varepsilon}=-\partial_t^{k+1}\psi_{\rm lin}^{\varepsilon}+\gamma\partial_x^{-1}\partial_t^{k+1}\eta_{\rm lin}^{\varepsilon}\in H^2.
\]
Thus $\partial_t^k \eta_{\rm lin}^{\varepsilon}\in H^4_0$, and iteration of these two relations gives smoothness in space. The smooth case therefore implies $\eta_{\rm lin}^{\varepsilon}=\psi_{\rm lin}^{\varepsilon}=0$ on $(\varepsilon,T-\varepsilon)$. Letting $\varepsilon\to0$ on compact subintervals of $(0,T)$ proves the result.
\end{proof}

If $\chi \partial_t\eta_{\rm lin}=0$ in distributions, Proposition~\ref{Prop-linear-unique} applies on any nonempty open interval where $\chi>0$.

\subsection{Observability and  energy decay}

We first estimate $N(\eta)\psi$ in $H^{-2}$ for the passage to the linearized limit. The Rellich identity gives the following bound in terms of the kinetic energy.

\begin{lemma}\label{Lemma-negative-N}
Under the assumptions of Lemma~{\rm\ref{Lemma-elliptic}}, we have
\begin{equation}\label{negative-N}
\|N(\eta)\psi\|_{H^{-2}}\le CK.
\end{equation}
\end{lemma}

\begin{proof}
Applying \eqref{rellich} and \eqref{bottomtrace} to a smooth periodic test function $\mu$, and using $2|\Phi_x\Phi_y|\le|\nabla\Phi|^2$, we obtain
\[
\left|\int_{-\pi}^{\pi}\mu N(\eta)\psi\d x\right|\le C\bigl(\|\mu\|_{L^\infty}+\|\mu_x\|_{L^\infty}\bigr)K\le C\|\mu\|_{H^2}K,
\]
where the last inequality follows from $H^2\hookrightarrow W^{1,\infty}$.
Taking the supremum over such $\mu$ with $\|\mu\|_{H^2}\le1$ proves
\eqref{negative-N}.
\end{proof}

For $s>5/2$, set
\[
X^s=H^{s+1/2}_0\times H^s_0,\quad\|(\eta,\widetilde\psi)\|_{X^s}=\|\eta\|_{H^{s+1/2}}+\|\widetilde\psi\|_{H^s}.
\]
This is the scale $H^{r+1/4}_0\times\dot H^{r-1/4}$ used in \cite{BMM24}, with $r=s+1/4$. Smallness in $X^s$ implies the geometric bounds required by Proposition~\ref{Prop-weak}. We use \eqref{weak-estimate} and Proposition~\ref{Prop-linear-unique} to prove the following observation inequality.

\begin{prop}\label{Prop-observation}
There exist $T_*>0$, $\varepsilon_*>0$ and $C_{\rm obs}\ge1$ such that every regular solution of \eqref{Waterwave} with feedback \eqref{Feedback-law} on $[0,T_*]$ for which
\begin{equation*}
\sup_{0\le t\le T_*}\|(\eta(t),\widetilde\psi(t))\|_{X^s}\le\varepsilon_*
\end{equation*}
satisfies
\begin{equation}\label{observability}
\mathcal E_\gamma(0)\le C_{\rm obs}\int_0^{T_*}\mathcal D(t)\d t.
\end{equation}
The constants depend only on $g,\kappa,h,\gamma,s$ and the cut-off.
\end{prop}

\begin{proof}
Choose $T_*>4C_w$, where $C_w$ is the constant in \eqref{weak-estimate}. If no pair $(\varepsilon_*,C_{\rm obs})$ has the required property for this $T_*$, we can choose regular solutions $(\eta_n,\psi_n)$ such that
\begin{equation}\label{contradiction-sequence}
\sup_{0\le t\le T_*}\|(\eta_n(t),\widetilde\psi_n(t))\|_{X^s}\rightarrow0,\quad\int_0^{T_*}\mathcal D_n(t)\d t<\frac1n\mathcal E_{\gamma,n}(0).
\end{equation}
For sufficiently large $n$, \eqref{coercivity} and \eqref{weak-estimate} apply. We normalize by setting
\[
r_n=\mathcal E_{\gamma,n}(0)^{1/2}>0,\quad
\widehat\eta_n=\frac{\eta_n}{r_n},\quad
\widehat\psi_n=\frac{\widetilde\psi_n}{r_n}.
\]
By \eqref{coercivity} and \eqref{contradiction-sequence}, $r_n\to0$, while the energy balance gives
\begin{equation}\label{normalized-energy}
1-\frac1n\le\frac{\mathcal E_{\gamma,n}(t)}{r_n^2}\le1\quad(0\le t\le T_*).
\end{equation}
Applying \eqref{coercivity} once more, we obtain
\begin{equation}\label{normalized-bound}
\sup_n\sup_{0\le t\le T_*}\left(\|\widehat\eta_n(t)\|_{H^1}+\|\widehat\psi_n(t)\|_{H^{1/2}}\right)<\infty.
\end{equation}
The feedback law also yields
\begin{equation}\label{pressure-limit}
\|P_{{\rm ext},n}/r_n\|_{L^2(0,T_*;L^2)}^2\le\|\chi\|_{L^\infty}\int_0^{T_*}\frac{\mathcal D_n(t)}{r_n^2}\d t\rightarrow 0.
\end{equation}

We next estimate the nonlinear terms in the normalized equations. By \eqref{DNcontinuity},
\begin{equation}\label{DN-limit}
\|(G(\eta_n)-G_0)\widehat\psi_n\|_{H^{-1/2}}\le C\|\eta_n\|_{W^{1,\infty}}\|\widehat\psi_n\|_{H^{1/2}}\rightarrow 0.
\end{equation}
The local term in the kinematic equation satisfies
\begin{equation}\label{kinematic-limit}
\frac{\eta_n\eta_{n,x}}{r_n}=\eta_n \partial_x\widehat\eta_n\rightarrow0\quad\text{in }L^\infty(0,T_*;L^2).
\end{equation}
For the corresponding term in the second equation, integration by parts gives
\[
\left|\int_{-\pi}^{\pi}\mu\eta_n \partial_x\widehat\psi_n\d x\right|=\left|\int_{-\pi}^{\pi}\widehat\psi_n\partial_x(\mu\eta_n)\d x\right|\le C\|\eta_n\|_{W^{1,\infty}}\|\widehat\psi_n\|_{L^2}\|\mu\|_{H^1}
\]
for every smooth periodic test function $\mu$. Hence
\begin{equation}
\frac{\eta_n\psi_{n,x}}{r_n}=\eta_n \partial_x\widehat\psi_n\rightarrow 0\quad\text{in }L^\infty(0,T_*;H^{-1}).
\end{equation}
For the curvature term, we write
\[
\frac{H(\eta_n)-\eta_{n,xx}}{r_n}=\partial_x\left[ (\partial_x\widehat\eta_n)\left((1+\eta_{n,x}^2)^{-1/2}-1\right)\right].
\]
Using $|(1+z^2)^{-1/2}-1|\le Cz^2$ for small $z$, we find
\begin{equation}
\left\|\frac{H(\eta_n)-\eta_{n,xx}}{r_n}\right\|_{H^{-1}}\le C\|\eta_{n,x}\|_{L^\infty}^2\|\partial_x\widehat\eta_n\|_{L^2}\rightarrow 0.
\end{equation}
Finally, $K_n\le Cr_n^2$ by \eqref{normalized-energy} and \eqref{coercivity}, so Lemma~\ref{Lemma-negative-N} gives
\begin{equation}\label{N-limit}
\left\|\frac{N(\eta_n)\psi_n}{r_n}\right\|_{H^{-2}}\le\frac{CK_n}{r_n}\le Cr_n\rightarrow0.
\end{equation}
The convergences in \eqref{DN-limit}--\eqref{N-limit} are uniform in time.

Subtracting the spatial mean from the second equation in \eqref{Waterwave} introduces the scalar
\[
c_n(t)=\left\langle N(\eta_n)\psi_n-\gamma\eta_n\psi_{n,x}+P_{{\rm ext},n}\right\rangle.
\]
By \eqref{mean-N}, \eqref{bottomtrace} and periodic integration by parts,
\[
\begin{aligned}
\frac{|\langle N(\eta_n)\psi_n\rangle|}{r_n}&\le\frac{CK_n}{r_n}\le Cr_n,\\
\frac{|\langle\eta_n\psi_{n,x}\rangle|}{r_n}
&= r_n|\langle (\partial_x\widehat\eta_n)\widehat\psi_n\rangle|\le Cr_n\|\partial_x\widehat\eta_n\|_{L^2}\|\widehat\psi_n\|_{L^2}.
\end{aligned}
\]
Together with \eqref{normalized-bound} and \eqref{pressure-limit}, these estimates give $c_n/r_n\to0$ in $L^2(0,T_*)$.
The normalized equations are
\[
\begin{aligned}
\partial_t\widehat\eta_n={}&G(\eta_n)\widehat\psi_n+\gamma\eta_n \partial_x\widehat\eta_n,\\
\partial_t\widehat\psi_n={}&-g\widehat\eta_n+\kappa \partial_x^2\widehat\eta_n+\gamma\partial_x^{-1}G(\eta_n)\widehat\psi_n-\frac{N(\eta_n)\psi_n}{r_n}\\
&+\kappa\frac{H(\eta_n)-\eta_{n,xx}}{r_n}+\gamma\eta_n \partial_x\widehat\psi_n-\frac{P_{{\rm ext},n}}{r_n}+\frac{c_n}{r_n}.
\end{aligned}
\]

The preceding estimates bound $\partial_t\widehat\eta_n,\partial_t\widehat\psi_n$ in $L^2(0,T_*;H^{-2})$. By \eqref{normalized-bound}, a subsequence converges weakly-* to $(\eta_{\rm lin},\psi_{\rm lin})$ in $L^\infty(0,T_*;H^1\times H^{1/2})$. For each fixed Fourier mode, the coefficients are bounded in $H^1(0,T_*)$ and therefore form a relatively compact set in $L^2(0,T_*)$. Moreover, if $\Pi_M$ denotes projection onto $|j|\le M$, \eqref{normalized-bound} gives
\[
\int_0^{T_*}\left(\|(1-\Pi_M)\widehat\eta_n\|_{L^2}^2+\|(1-\Pi_M)\widehat\psi_n\|_{L^2}^2\right)\d t\le CT_*(M^{-2}+M^{-1}).
\]
Combining the compactness of the finite Fourier sums with this tail estimate, we obtain, after passing to a further subsequence,
\begin{equation}\label{compact-convergence}
(\widehat\eta_n,\widehat\psi_n)\rightarrow(\eta_{\rm lin},\psi_{\rm lin}) \quad\text{strongly in }L^2(0,T_*;L^2\times L^2).
\end{equation}

The nonlocal term converges by \eqref{DN-limit} and \eqref{compact-convergence}, since the maps
$\partial_x^{-1}:H^{-1/2}_0\rightarrow H^{1/2}_0$ and $\partial_x^{-1}G_0:L^2_0\rightarrow L^2_0$ are bounded. 
Using also \eqref{kinematic-limit}--\eqref{N-limit}, \eqref{pressure-limit} and $c_n/r_n\to0$, we pass to the limit in the normalized equations and conclude that $(\eta_{\rm lin},\psi_{\rm lin})$ solves \eqref{linear-system} in distributions, with both components of zero mean. The feedback law gives
\[
\chi\partial_t\widehat\eta_n=P_{{\rm ext},n}/r_n\rightarrow0\quad\text{in }L^2(0,T_*;L^2).
\]
Thus $\chi\partial_t\eta_{\rm lin}=0$, and Proposition~\ref{Prop-linear-unique} yields $\eta_{\rm lin}=\psi_{\rm lin}=0$. By \eqref{compact-convergence},
\begin{equation}\label{compact-remainder-zero}
\int_0^{T_*}\|\widehat\psi_n(t)\|_{L^2}^2\d t\rightarrow0.
\end{equation}
Applying \eqref{weak-estimate} to $(\eta_n,\psi_n)$ and dividing by $r_n^2$, we obtain from \eqref{normalized-energy} and \eqref{contradiction-sequence}
\[
\left(1-\frac1n\right)T_*\le C_w\left[2+\frac1n+\int_0^{T_*}\|\widehat\psi_n(t)\|_{L^2}^2\d t\right].
\]
Letting $n\to\infty$ and using \eqref{compact-remainder-zero} gives $T_*\le2C_w$, contrary to our choice of $T_*$. This proves \eqref{observability}.
\end{proof}

Applying the observation inequality on successive time intervals gives the following bound for the time integral of the energy.

\begin{lemma}\label{Lemma-integrated}
There is $C>0$ such that every regular solution of \eqref{Waterwave} with feedback \eqref{Feedback-law} on $[0,T]$ for which
\[
\sup_{0\le t\le T}\|(\eta(t),\widetilde\psi(t))\|_{X^s}\le\varepsilon_*
\]
satisfies
\begin{equation}\label{integrated-energy}
\int_0^T\mathcal E_\gamma(t)\d t\le C\mathcal E_\gamma(0).
\end{equation}
The constant $C$ is independent of $T$ and of the solution.
\end{lemma}

\begin{proof}
Let $N=\lfloor T/T_*\rfloor$. For $N\ge1$, we apply \eqref{observability} by time translation on each interval $[kT_*,(k+1)T_*]$, $k=0,\ldots,N-1$. The monotonicity of the energy and \eqref{dissipation} then give
\[
\int_{kT_*}^{(k+1)T_*}\mathcal E_\gamma(t)\d t\le T_*\mathcal E_\gamma(kT_*)\le T_*C_{\rm obs}
\bigl(\mathcal E_\gamma(kT_*)-\mathcal E_\gamma((k+1)T_*)\bigr).
\]
Summing and estimating the remaining interval by monotonicity, we obtain
\[
\int_0^T\mathcal E_\gamma(t)\d t\le T_*C_{\rm obs} \bigl(\mathcal E_\gamma(0)-\mathcal E_\gamma(NT_*)\bigr)+T_*\mathcal E_\gamma(NT_*)\le T_*C_{\rm obs}\mathcal E_\gamma(0),
\]
since $C_{\rm obs}\ge1$. For $T<T_*$, the same bound follows directly from monotonicity. Thus \eqref{integrated-energy} holds with $C=T_*C_{\rm obs}$.
\end{proof}

We now state and prove the main result.

\begin{theorem}\label{Thm-stabilization}
Let $g,\kappa,h>0$, $\gamma\in\R$ and $s>5/2$ be fixed. Choose $m,\chi$ as in \eqref{cutoff}, and assume \eqref{capillary-cutoff-condition}. There exist $\varepsilon_*>0$ and $C>0$ such that the following holds. Let $(\eta,\psi,P_{\rm ext})$ be a regular solution of \eqref{Waterwave} with feedback \eqref{Feedback-law} on $[0,T]$, and suppose that
\begin{equation}
\sup_{0\le t\le T}\left(\|\eta(t)\|_{H^{s+1/2}}+\|\psi(t)-\langle\psi(t)\rangle\|_{H^s}\right)\le\varepsilon_*.
\end{equation}
Then
\begin{equation}\label{main-integral}
\int_0^T\mathcal E_\gamma(t)\d t\le C\mathcal E_\gamma(0),\quad\mathcal E_\gamma(T)\le\frac CT\mathcal E_\gamma(0) \quad(T>0).
\end{equation}
Moreover,
\begin{equation}\label{main-exponential}
\mathcal E_\gamma(t) \le 2\exp\left(-\frac{\log2}{2C}t\right)\mathcal E_\gamma(0),\quad 0\le t\le T.
\end{equation}
The constants depend only on $g,\kappa,h,\gamma,s$ and the chosen cut-off, and are independent of $T$ and of the solution.
\end{theorem}

\begin{proof}
Choose $\varepsilon_*$ as in Proposition~\ref{Prop-observation}, small enough for the preceding geometric estimates to apply. Lemma~\ref{Lemma-integrated} gives the first estimate in \eqref{main-integral}, while monotonicity yields
\[
T\mathcal E_\gamma(T)\le\int_0^T\mathcal E_\gamma(t)\d t\le C\mathcal E_\gamma(0),
\]
which proves the second estimate.
Set $T_0=2C$. By time translation, the second estimate in \eqref{main-integral} applies on each interval $[a,a+T_0]\subset[0,T]$ and gives
\[
\mathcal E_\gamma(a+T_0)\le\frac{C}{T_0}\mathcal E_\gamma(a)=\frac12\mathcal E_\gamma(a).
\]
For $t\in[0,T]$, take $n=\lfloor t/T_0\rfloor$. Iteration and
monotonicity imply
\[
\mathcal E_\gamma(t)\le\mathcal E_\gamma(nT_0)\le 2^{-n}\mathcal E_\gamma(0)\le 2\exp\left(-\frac{\log2}{T_0}t\right)\mathcal E_\gamma(0),
\]
since $n\ge t/T_0-1$. This is \eqref{main-exponential}.
\end{proof}

\begin{remark}
Thanks to \eqref{coercivity}, i.e. the equivalence of $H^1\times H^{1/2}$ norm of $(\eta,\tilde\psi)$ and energy $\mathcal E_\gamma$, one can further derive an exponential decay for the former. More precisely, there is $C_1>0$ such that, with $\alpha=\log2/(2C)$,
\[
\|\eta(t)\|_{H^1}+\|\widetilde\psi(t)\|_{H^{1/2}}\le C_1e^{-\alpha t/2}\left(\|\eta(0)\|_{H^1}+\|\widetilde\psi(0)\|_{H^{1/2}}\right).
\]
In addition,  we  have
\[
\|\nabla\Phi(t)\|_{L^2(D_\eta(t))}\le C_1e^{-\alpha t/2}\left(\|\eta(0)\|_{H^1}+\|\widetilde\psi(0)\|_{H^{1/2}}\right),
\]
since $
\|\nabla\Phi(t)\|_{L^2(D_\eta(t))}^2=2K(t)\le C\mathcal E_\gamma(t).
$
\end{remark}

\medskip

\subsection*{Acknowledgements}
The authors would like to thank Professor Jifeng Chu for valuable discussions and suggestions during the preparation of the paper.

Qixing Ding is partially supported by the National Natural Science Foundation of China (Grant No. 12601316) and the Shanghai Pujiang Program (Grant No. 25PJA044). Jia-Cheng Zhao is partially supported by the National Natural Science Foundation of China (Grant No. 12601877).

\subsection*{Use of AI declaration}

OpenAI’s ChatGPT was used as an aid in mathematical discussions, language polishing and proof checking. The authors conceived and developed the work and finalized the manuscript, and take full responsibility for the mathematical content.

\end{document}